\documentclass[10pt]{amsart}
\usepackage{amssymb}
\usepackage{enumerate}
\usepackage{epsf}
\usepackage{psfrag}
\DeclareGraphicsExtensions{.eps,.art,.ART,.ps}
\newcommand{\bbR}{\mathbb{R}}      

\newtheorem{Thm}{Theorem}[section]

\begin{document}
\begin{abstract}
In this note we derive a new Minkowski-type inequality for closed convex surfaces in the hyperbolic $3$-space. The inequality is obtained by explicitly computing the area of the family of surfaces obtained from the normal flow and then applying the isoperimetric inequality. Using the same method, we also we give elementary proofs of the classical Minkowski inequalities for closed convex surfaces in the Euclidean $3$-space and in the $3$-sphere.
\end{abstract}
%
%
\title[A Minkowski-type inequality for convex surfaces in $H^3$]{A Minkowski-type inequality for convex surfaces in the hyperbolic $3$-space}
\author{Jos\'{e} Nat\'{a}rio}
\address{Centro de An\'alise Matem\'atica, Geometria e Sistemas Din\^amicos, Departamento de Matem\'atica, Instituto Superior T\'ecnico, 1049-001 Lisboa, Portugal}
\thanks{Partially supported by FCT (Portugal).}
\maketitle
%
%
%
\section*{Introduction}
The classical Minkowski inequality for closed convex surfaces $S \subset \bbR^3$ (that is, smooth surfaces which are boundaries of convex bounded open sets) reads
\begin{equation} \label{Minkowski}
\int_{S} H \geq \sqrt{16 \pi |S|}
\end{equation}
with equality if and only if $S$ is a sphere, where $H \geq 0$ is the mean curvature\footnote{We define the mean curvature to be the trace of its second fundamental form. Notice that many authors define it to be {\em half} this trace, leading to a different numerical factor in the Minkowski inequality.} of $S$ and $|S|$ is its surface area. It was used by Minkowski himself to give a proof of the isoperimetric inequality for closed convex surfaces \cite{Minkowski03, Osserman78}, by Gibbons to prove the Penrose inequality for null dust shells in Minkowski spacetime \cite{Gibbons97, Mars09}, and by Lam to prove the Penrose inequality for graphs \cite{Lam10}. Minkowski-type inequalities for closed convex surfaces in hyperbolic space were proved by Gallego and Solanes \cite{GS05} and de Lima and Gir\~ao \cite{dLG12}, who then used them to prove the Penrose inequality for asymptotically hyperbolic graphs. Similar inequalities were obtained by Brendle, Hung and Wang for the the Anti-de Sitter-Schwarzschild manifold \cite{BHW12}.

In this note we reverse Minkowski's idea and prove the Minkowski inequality starting with the isoperimetric inequality. The idea is the following: using the Gauss-Bonnet theorem, we obtain (Section~\ref{section1}) an explicit expression for the area of the surfaces $S_t$ obtained from a compact co-orientable surface $S$ by flowing a distance $t$ along the orthogonal geodesics. If $S$ is convex and we take the outward co-orientation then this flow is well defined for all $t \geq 0$; we then observe (Section~\ref{section2}) that if the Minkowski inequality did not hold for $S$ then $S_t$ would violate the isoperimetric inequality for sufficiently large $t$ (note that we do not obtain the rigidity statement). 

Since the formulas obtained in Section~\ref{section1} work for surfaces in any constant curvature $3$-manifold, we can use the same method in the hyperbolic space. In Section~\ref{section3} we prove:

\begin{Thm}\label{Thm1}
Let $S$ be a closed convex surface in the hyperbolic $3$-space $H^3$ (that is, $S$ is a smooth boundary of a geodesically convex bounded open set). Then
\begin{equation} \label{MinkowskiHyp}
\int_{S} H \geq 4 V + 4\pi\log\left(1 + \frac{|S|}{2\pi} + \frac{1}{4\pi}\int_{S} H \right),
\end{equation}
where $H$ is the mean curvature of $S$, $|S|$ is its surface area and $V$ is the volume enclosed by $S$.
\end{Thm}

Although the equality holds when $S$ is a sphere, we do not obtain a rigidity statement (nor do we know whether it is true). We note that the weaker inequality $\int_{S} H \geq 4 V$ was proved in \cite{GS05} (cf.~footnote 1 when comparing). Still in Section~\ref{section3} we explain how \eqref{MinkowskiHyp}, which is seemingly quite different from the Minkowski inequality \eqref{Minkowski}, reduces to it for surfaces $S$ with small diameter (as compared to the radius of curvature of the hyperbolic space). 

For completeness, the formulas obtained in Section~\ref{section1} are used in Section~\ref{section4} to prove:

\begin{Thm}\label{Thm2}
Let $S$ be a closed convex surface in the $3$-sphere $S^3$ (that is, $S$ is a smooth boundary of a geodesically convex open set). Then
\begin{equation} \label{MinkowskiSph}
\int_{S} H \geq \sqrt{16 \pi |S| \left( 1 -  \frac{|S|}{4\pi} \right)},
\end{equation}
where $H$ is the mean curvature of $S$ and $|S|$ is its surface area. Moreover, equality holds if and only if $S$ is a sphere.
\end{Thm}

This result was stated without proof in \cite{Blaschke38}, and was first proved in \cite{Knothe52}. Our proof is essentially the one given in \cite{Santalo63}. In this case we do obtain a rigidity result; this is because for the $3$-sphere the potential violation of the isoperimetric inequality occurs at a finite distance from $S$, rather than at infinity, as is the case for the Euclidean or the hyperbolic spaces.

Finally, in Section~\ref{section5} we compare Theorems~\ref{Thm1} and \ref{Thm2} and the classical Minkowski inequality, and discuss other Minkowski-type inequalities for closed convex surfaces in $H^3$.
%
%
%
\section{Normal flow}\label{section1}
Let $(M,g)$ be a Riemannian $(n+1)$-dimensional manifold and $S \subset M$ a compact co-orientable hypersurface. Given a choice of a smooth unit normal $\nu$ on $S$ we define (for sufficiently small $\varepsilon > 0$) the map $\exp:(-\varepsilon, \varepsilon) \times S \to M$ by the formula $\exp(t,p)=c_p(t)$, where $c_p:(-\varepsilon, \varepsilon) \to M$ is the geodesic with initial condition $\dot{c}_p(0)=\nu_p$. We can choose $\varepsilon$ so that each of the maps $\exp_t:S \to M$ defined by $\exp_t(p)=\exp(t,p)$ is an embedding, and so each set $S_t = \exp_t(S)$ is a compact co-orientable hypersurface with unit normal given by $\nu_{\exp_t(p)}=\dot{c}_p(t)$. As is well known (see for instance \cite{CLN07}), the $n$-dimensional area of these hypersurfaces satisfies
\[
\frac{d}{dt} |S_t| = \int_{S_t} H
\]
and
\[
\frac{d^2}{dt^2} |S_t| = \int_{S_t} \left( Ric(\nu,\nu) - R + \bar{R} \right),
\]
where $H$ and $\bar{R}$ are the mean and scalar curvatures of $S_t$, and $Ric$ and $R$ are the Ricci and scalar curvatures of $M$. In the particular case when $M$ is a manifold of constant curvature $K$ we have $Ric(\nu,\nu)=nK$ and $R=(n+1)nK$, yielding
\[
\frac{d^2}{dt^2} |S_t| = -n^2 K |S_t| + \int_{S_t} \bar{R}.
\]
If $n=2$, we obtain from the Gauss-Bonnet Theorem
\[
\frac{d^2}{dt^2} |S_t| = -4 K |S_t| + 4\pi\chi(S_t),
\]
where $\chi(S_t)$ is the Euler characteristic of $S_t$. In the case when $S$ (and therefore $S_t$) is homeomorphic to $S^2$, this becomes
\begin{equation} \label{main}
\frac{d^2}{dt^2} |S_t| = -4 K |S_t| + 8\pi.
\end{equation}
%
%
%
\section{Euclidean $3$-space}\label{section2}
As a warm-up exercise, we give an elementary proof of the Minkowski inequality for closed convex surfaces in the Euclidean $3$-space $\bbR^3$ (that is, smooth surfaces which are boundaries of convex bounded open sets). Notice that such surfaces are automatically homeomorphic to $S^2$, and their normal flow is well defined for $t \geq 0$. If we set $A(t)=|S_t|$ then we have from \eqref{main} with $K=0$
\[
\ddot{A}(t)=8\pi,
\]
which can be immediately integrated to
\[
A(t)=4\pi t^2 + \dot{A}_0 t + A_0,
\]
where 
\[
A_0=A(0)=|S|
\]
and 
\[
\dot{A}_0=\dot{A}(0)=\int_S H.
\]
The volume $V(t)$ of the convex open set bounded by $S_t$ varies as
\[
\dot{V}(t)=A(t),
\]
and so
\[
V(t)=\frac43 \pi t^3 + \frac12 \dot{A}_0 t^2 + A_0 t + V_0,
\]
where $V_0=V(0)$. A straightforward computation yields
\[
A^3 - 36\pi V^2 =  3\pi ({\dot{A}_0}^2 - 16 \pi A_0 ) \, t^4 + P(t),
\]
where $P$ is a polynomial of degree $3$. Since the isoperimetric inequality
\[
A^3 \geq 36\pi V^2
\]
must hold for all $t \geq 0$, we conclude that
\[
{\dot{A}_0}^2 \geq 16 \pi A_0,
\]
which is Minkowski's inequality \eqref{Minkowski}.

Note however that this method does not yield the rigidity statement (that if the equality holds in \eqref{Minkowski} then $S$ is a sphere).
%
%
%
\section{Hyperbolic $3$-space}\label{section3}
We now turn to closed convex surfaces in the hyperbolic $3$-space $H^3$ (that is, smooth surfaces which are boundaries of geodesically convex bounded open sets). Notice that such surfaces are automatically homeomorphic to $S^2$, and their mean curvature is nonnegative. Since $H^3$ has negative curvature $K=-1$, the normal flow cannot develop conjugate points, and a simple argument involving the sum of the internal angles of a geodesic triangle shows that different normal geodesics cannot cross. We conclude that the normal flow is is well defined for $t \geq 0$. We now have from \eqref{main}, again setting $A(t)=|S_t|$,
\[
\ddot{A}(t)=4A(t)+8\pi,
\]
which can be immediately integrated to
\[
A(t)=2\pi R e^{2t} + 2 \pi T e^{-2t} - 2\pi,
\]
where
\[
\begin{cases}
\displaystyle R + T = 1 + \frac{A_0}{2\pi} \\
\displaystyle R - T = \frac{\dot{A}_0}{4\pi}
\end{cases}
\]
($R \geq 1$, $T \leq R$), and again
\[
A_0=A(0)=|S|
\]
and 
\[
\dot{A}_0=\dot{A}(0)=\int_S H.
\]
The volume $V(t)$ of the convex open set bounded by $S_t$ varies as
\[
\dot{V}(t)=A(t),
\]
and so
\[
V(t)= \pi R e^{2t} - \pi T e^{-2t} - 2\pi t + \pi (T-R) + V_0,
\]
where $V_0=V(0)$.

The formulas for the area and volume of a sphere of radius $r$ can be obtained from the limit when $S$ approaches a single point ($A_0=\dot{A}_0=V_0=0$):
\[
A(r)=2\pi \cosh(2r) - 2\pi
\]
and
\[
V(r)=\pi \sinh(2r) - 2\pi r.
\]
The radius $r(t)$ of the sphere with the same area as $S_t$ can then be obtained by solving
\[
\cosh(2r) = R e^{2t} + T e^{-2t}.
\]
As $t \to + \infty$ we have
\[
2r \sim \log (2R) + 2t
\]
with exponentially small error; therefore
\begin{align*}
V(t) & \sim \pi \sinh(2r) - 2\pi r + \pi \log (2R) + \pi (T-R) + V_0 \\
& \sim V(r) + \pi\log \left(1 + \frac{A_0}{2\pi} + \frac{\dot{A}_0}{4\pi}\right) - \frac{\dot{A}_0}{4} + V_0,
\end{align*}
again with exponentially small error. Since by the isoperimetric inequality we must have $V(t) \leq V(r(t))$ for all $t \geq 0$, we conclude that
\[
\frac{\dot{A}_0}{4\pi} \geq \frac{V_0}{\pi} + \log \left(1 + \frac{A_0}{2\pi} + \frac{\dot{A}_0}{4\pi}\right),
\]
which is exactly inequality \eqref{MinkowskiHyp} in Theorem~\ref{Thm1}. 

It is interesting to see how this inequality reduces to the Euclidean Minkowski inequality \eqref{Minkowski} for small surfaces. One first notices that it implies
\[
\frac{\dot{A}_0}{4\pi} \geq  \log \left(1 + \frac{A_0}{2\pi} + \frac{\dot{A}_0}{4\pi}\right) \Leftrightarrow \exp\left(\frac{\dot{A}_0}{4\pi}\right) - \frac{\dot{A}_0}{4\pi} - 1 \geq  \frac{A_0}{2\pi}.
\]
Noting that
\[
e^x - x - 1 = \frac{x^2}{2} + \frac{x^3}{6} + \ldots,
\]
we see that to second order in $\dot{A}_0$ the inequality is just
\[
{\dot{A}_0}^2 \geq 16 \pi A_0,
\]
which is the Euclidean Minkowski inequality \eqref{Minkowski}. Note that $\dot{A}_0$ is of the order of the diameter of the set bounded by $S$, and so requiring $\dot{A}_0 \ll 1$ (so that the second order approximation is accurate) is the same as requiring this diameter to be much smaller then the radius of curvature of the hyperbolic space.

A straightforward calculation shows that if $S$ is a sphere then the equality in \eqref{MinkowskiHyp} holds. However, it is not clear if one can make a rigidity statement (that if the equality holds in \eqref{MinkowskiHyp} then $S$ is a sphere).
%
%
%
\section{$3$-sphere}\label{section4}
Finally, we turn to closed convex surfaces in the $3$-space $S^3$ (that is, smooth surfaces which are boundaries of geodesically convex open sets\footnote{Recall that a subset $S \subset M$ of a Riemannian manifold $(M,g)$ is geodesically convex if for any two points $p, q \in S$ there exists a {\em minimizing} geodesic in S joining p to q.}). Notice that these surfaces are automatically homeomorphic to $S^2$, and their mean curvature is nonnegative. Since $S^3$ has constant curvature $K=1$, any conjugate points of the normal flow lie at a distance greater than or equal to $\frac{\pi}2$ from $S$, and a simple argument involving the sum of the internal angles of a geodesic triangle contained in a hemisphere shows that different normal geodesics can only cross at a distance greater than or equal to $\frac{\pi}2$ from $S$ along one of them. We conclude that the normal flow is is well defined for $0 \leq t < \frac{\pi}{2}$. We now have from \eqref{main}, again setting $A(t)=|S_t|$,
\[
\ddot{A}(t)=-4A(t)+8\pi,
\]
which can be immediately integrated to
\[
A(t)=2\pi - 2\pi R \cos(2t + \theta),
\]
where
\[
\begin{cases}
\displaystyle R \cos \theta = 1 - \frac{A_0}{2\pi} \\
\displaystyle R \sin \theta = \frac{\dot{A}_0}{4\pi}
\end{cases}
\]
($R\geq 0$, $0 \leq \theta \leq \pi$), and again
\[
A_0=A(0)=|S|
\]
and 
\[
\dot{A}_0=\dot{A}(0)=\int_S H.
\]
The volume $V(t)$ of the convex open set bounded by $S_t$ varies as
\[
\dot{V}(t)=A(t),
\]
and so
\[
V(t)=  2\pi t - \pi R \sin(2t + \theta) + \pi R \sin\theta + V_0,
\]
where $V_0=V(0)$.

The formulas for the area and volume of a sphere of radius $r$ can be obtained from the limit when $S$ approaches a single point ($A_0=\dot{A}_0=V_0=0$):
\[
A(r)=2\pi - 2\pi \cos(2r)
\]
and
\[
V(r)=2\pi r - \pi \sin(2r).
\]

The function $A(t)$ reaches its maximum value $2\pi(R+1)$ for $t=(\pi-\theta)/2$. The function $V(t)$, on the other hand, is increasing and
\[
V\left(\frac{\pi}2\right) =  \pi^2 + 2\pi R \sin\theta + V_0 > \pi^2.
\]
So if $R<1$ there would exist a surface with area smaller that $4\pi$ bounding a volume equal to $\pi^2$, violating the isoperimetric inequality. Therefore we must have
\[
R \geq 1 \Leftrightarrow  \left(1 - \frac{A_0}{2\pi}\right)^2 + \left(\frac{\dot{A}_0}{4\pi}\right)^2 \geq 1,
\]
and so
\[
{\dot{A}_0}^2 \geq 16 \pi A_0 \left( 1 -  \frac{A_0}{4\pi} \right),
\]
which is exactly inequality \eqref{MinkowskiSph} in Theorem~\ref{Thm2}. If $R=1$, on the other hand, we will have $V(t)=\pi^2$ with $A(t) \leq 4\pi$, and so in fact we must have $A(t)=4\pi$ and $S_t$ must be a unit $2$-sphere, meaning that $S$ is also a sphere. This yields the rigidity statement in Theorem~\ref{Thm2}.
%
%
%
\section{Other inequalities}\label{section5}
If we define $r(t)$ as the radius of the sphere with the same area as $S_t$, a straightforward calculation shows that inequalities \eqref{Minkowski} and \eqref{MinkowskiSph} for closed convex surfaces in $\bbR^3$ and $S^3$ can be reinterpreted as saying that
\begin{equation} \label{r}
\frac{dr}{dt}(0) \geq 1.
\end{equation}
For the hyperbolic space $H^3$, \eqref{r} yields the tempting inequality
\begin{equation}\label{false}
{\dot{A}_0}^2 \geq 16 \pi A_0 \left( 1 +  \frac{A_0}{4\pi} \right),
\end{equation}
which unfortunately cannot be true. This can be seen by considering a convex surface approximating a geodesic disk of radius $R$, for which $A_0$ approaches $4\pi(\cosh R-1)$ and $\dot{A}_0$ approaches $2 \pi^2 \sinh R$, thus violating \eqref{false} for sufficiently large $R$. This observation is due to Naveira and Solanes \cite{NS09}.

Although \eqref{false} does not hold, the weaker inequality
\[
{\dot{A}_0}^2 \geq {A_0}^2,
\]
does, as shown in \cite{GS05} (cf.~footnote 1 when comparing). Another weaker inequality suggested by \eqref{false} is just the classical Minkowski inequality
\[
{\dot{A}_0}^2 \geq 16 \pi A_0.
\]
It is not clear whether it applies to the hyperbolic space $H^3$ as well.
%
%
%
\section*{Acknowledgments}
I thank Gil Solanes for bringing to my attention important references, and Levi de Lima for clarifying some issues related to his work.
%
%
%

\end{document}